\documentclass[10pt]{amsart}
\usepackage{mathrsfs}
\usepackage{amsfonts} 
\usepackage{graphicx,latexsym,bm,amsmath,amssymb,verbatim,multicol,lscape}
\theoremstyle{definition}

\theoremstyle{remark}

\numberwithin{equation}{section}

\begin{document}
\title[New results on permutation polynomials over finite fields]
{New results on permutation polynomials over finite fields}%
\author{Xiaoer Qin}
\address{Mathematical College, Sichuan University, Chengdu 610064, P.R. China;
and College of Mathematics and Computer Science, Yangtze Normal University,
Chongqing 408100, P.R. China}
\email{qincn328@sina.com}
\author{Guoyou Qian}
\address{Mathematical College, Sichuan University, Chengdu 610064, P.R. China}
\email{qiangy1230@163.com, qiangy1230@gmail.com}
\author{Shaofang Hong$^*$}
\address{Mathematical College, Sichuan University, Chengdu 610064, P.R. China}
\email{sfhong@scu.edu.cn, s-f.hong@tom.com, hongsf02@yahoo.com}
\thanks{$^*$S. Hong is the corresponding author and was supported partially
by National Science Foundation of China Grant \#11371260.}
\date{\today}%
\begin{abstract}
In this paper, we get several new results on permutation polynomials
over finite fields. First, by using the linear translator, we
construct permutation polynomials of the forms $L(x)+\sum_{j=1}^k
\gamma_jh_j(f_j(x))$ and $x+\sum_{j=1}^k\gamma_jf_j(x)$. These
generalize the results obtained by Kyureghyan in 2011. Consequently, we
characterize permutation polynomials of the form $L(x)+\sum_{i=1}
^l\gamma_i {\rm Tr}_{{\bf F}_{q^m}/{\bf F}_{q}}(h_i(x))$, which
extends a theorem of Charpin and Kyureghyan obtained in 2009.
\end{abstract}

\maketitle

\section{\bf Introduction}
Let $p$ be a prime and $q=p^{n}$ for some $n\in \textbf{Z}^+$ (the set of positive integers).
Let ${\bf F}_p$ be the prime field and ${\bf F}_{q}$ denote the finite field with $q$ elements.
Throughout ${\bf F}_{q}^{*}:={\bf F}_{q}\setminus \{0\}$ and ${\bf F}_{q}[x]$ represents the
ring of polynomials over ${\bf F}_{q}$ in the indeterminate $x$. A polynomial $f(x)\in {\bf
F}_{q}[x]$ is called a {\it permutation polynomial} of ${\bf F}_{q}$ if $f(x)$ induces
a permutation of ${\bf F}_{q}$. More information of permutation polynomials can be
found in the book of Lidl and Niederreiter \cite{[LN]}. Permutation polynomials have
many important applications in coding theory \cite{[L]}, cryptography \cite{[SH]}
and combinatorial design theory. The problem of constructing new classes of permutation
polynomials over finite fields has generated much interest, see the open problems
in \cite{[LM]}. Wan and Lidl \cite{[WL]}, Masuda and Zieve \cite{[MZ]} and Zieve \cite{[Z2]}
constructed permutation polynomials of the form $x^{r}f(x^{(q-1)/d})$ and studied
their group structure. Zieve \cite{[Z]} characterized the permutation
polynomial of the form $x^r(1+x^v+x^{2v}+...+x^{kv})^t$. Recently, by
using a powerful lemma, Zieve \cite{[Z3], [Z4]} got some new permutation
polynomials over finite fields. Ayad,
Belghaba and Kihel \cite{[ABK]} obtained some permutation binomials
and proved the bound of $p$, if $ax^n+x^m$ permutes ${\bf F}_{p}$.
Hou \cite{[H]} characterized two new classes of permutation polynomials over finite fields.

Let $m>1$ be a given integer. Throughout $\textup{Tr}_{{\bf F}_{q^m}/{\bf F}_{q}}(x)$
denotes  the {\it trace} from ${\bf F}_{q^m}$ to ${\bf F}_{q}$, that is,
$\textup{Tr}_{{\bf F}_{q^m}/{\bf F}_{q}}(x)=x+x^{q}+\cdots+x^{q^{m-1}}.$ In particular,
one has $\textup{Tr}_{{\bf F}_{q}/{\bf F}_{p}}(x)=x+x^{p}+\cdots+x^{p^{n-1}}.$
A polynomial of the form $L(x)=\sum^{m-1}_{i=0}a_{i}x^{q^{i}}\in {\bf F}_{q^m}[x]$ is
called a \emph{linearized polynomial} over ${\bf F}_{q^m}$. It is
well known that a linearized polynomial $L(x)$ is a permutation
polynomial of ${\bf F}_{q^m}$ if and only if the set of roots of $L(x)$
in ${\bf F}_{q^m}$ equals $\{0\}$ (see, for example, Theorem 7.9 of \cite{[LN]}).
Using the trace function and linearized polynomials, a number of classes
of permutation polynomials were constructed. Qin and Hong \cite{[QH]}
constructed permutation polynomials of the form $\sum_{i=1}^k(L_{i}(x)
+\gamma_i)h_i(B(x))$, where $L_i(x)$ and $B(x)$ are linearized polynomials.
Marcos \cite{[M]} obtained permutation polynomials of the form
$L(x)+\gamma h(\textup{Tr}_{{\bf F}_{q^m}/{\bf F}_{q}}(x))$. Zieve
\cite{[Z1]} presented rather more general versions of the first four
constructions from \cite{[M]}.

The linear translator is a
powerful technique to construct permutation polynomials. There are several classes of
permutation polynomials constructed by the linear translator. Charpin
and Kyureghyan \cite{[CK2]} studied permutation polynomials of the shape
$G(x)+\gamma \textup{Tr}_{{\bf F}_{2^n}/{\bf F}_{2}}(H(x))$ over ${\bf F}_{2^n}$.
Using the functions having linear translators, Charpin
and Kyureghyan \cite{[CK]} introduced an effective method to construct
permutation polynomials of the shape $G(x)+\gamma \textup{Tr}_{{\bf F}_{q}
/{\bf F}_{p}}(H(x))$ over ${\bf F}_{q}$, where $G(x)$ is either a permutation
or linearized polynomial. In \cite{[K]}, Kyureghyan further constructed
permutation polynomials of the forms $x+\gamma f(x)$ and $L(x)+\gamma h(f(x))$,
where $f(x)$ has a linear translator. Using linear translators, Qin and
Hong \cite{[QH]} characterized a class of permutation polynomials of
the form  $L_1(x)+L_{2}(\gamma)h(f(x))$, which generalizes a result of \cite{[K]}.

In this paper, our main goal is to construct some new permutation polynomials
over finite fields. First, in Section 2, by using the linear
translator, we characterize permutation polynomials of the forms
$L(x)+\sum_{j=1}^k\gamma_jh_j(f_j(x))$ and $x+\sum_{j=1}^k\gamma_jf_j(x)$. These
generalize the theorems of Kyureghyan \cite{[K]} obtained in 2011. Consequently, in
Section 3, we characterize permutation polynomials of the shape
$L(x)+\sum_{i=1}^l\gamma_i \textup{Tr}_{{\bf F}_{q^m}/{\bf F}_{q}}(h_i(x))$.
This extends a result due to Charpin and Kyureghyan \cite{[CK]}.

\section{\bf Constructing permutation polynomials  by linear translators}

In this section, we use the linear translator to construct two new classes of
permutation polynomials over finite fields. We first recall the definition of
linear translator as follows: \\

\noindent{\bf Definition 2.1.} \cite{[K]} Let $f: {\bf F}_{q^m}\to
{\bf F}_{q}$, $a\in {\bf F}_{q}$ and $\alpha $ be a nonzero element
in ${\bf F}_{q^m}$. If $f(x+u\alpha)-f(x)=ua$ for all $x\in {\bf
F}_{q^m}$ and all $u\in {\bf F}_{q}$, then we say that $\alpha $ is an
{\it $a$-linear translator} of the function $f$. In particular,
$a=f(\alpha)-f(0).$\\

In \cite{[CK]}, the functions holding a linear translator are
characterized as follows:\\

\noindent{\bf Lemma 2.1.} \cite{[CK]} {\it A mapping $f:{\bf
F}_{q^m}\to {\bf F}_{q}$ has a linear translator if and only if
there is a non-bijective linearized polynomial $L(x)\in {\bf
F}_{q^m}[x]$ such that $f(x)=\textup{Tr}_{{\bf F}_{q^m}/{\bf
F}_{q}}(\beta x+H(L(x)))$ for some mapping $H:{\bf F}_{q^m}\to {\bf
F}_{q^m}$ and $\beta\in {\bf F}_{q^m}$.}\\

Ling and Qu \cite {[LQ]} answered an open problem of \cite{[CK]}
and present a method to construct explicitly linearized polynomials with
kernel of any given dimension. We can now use the linear
translator to construct permutation polynomials and
give the first main result of this paper as follows.\\

\noindent{\bf Theorem 2.1.} {\it Let $k$ be a positive integer.
Let $L:{\bf F}_{q^m}\to {\bf F}_{q^m}$ be a linearized polynomial
such that ${\rm dim (Ker}(L))=k$ and
${\rm Ker}(L)\cap {\rm Im}(L)=\{0\}$. Let $\{\gamma_1,...,\gamma_k\}$ be a
basis of ${\rm Ker}(L)$ over ${\bf F}_q$ and $h_1(x), ..., h_k(x)\in {\bf F}_{q}[x]$
be permutation polynomials of ${\bf F}_{q}$. For any integers $i$ and $j$ with
$1\leq i,j\leq k$, let $b_{ij}\in {\bf F}_q$ and $\gamma_i$ be a
$b_{ij}$-linear translator of $f_j:{\bf F}_{q^m}\to {\bf F}_{q}$.
Then $F(x):=L(x)+\sum_{j=1}^k\gamma_jh_j(f_j(x))$ is a
permutation polynomial of ${\bf F}_{q^m}$ if and only if $\det\big(b_{ij}\big)_
{1\leq i, j\leq k}\neq0$.}
\begin{proof}
First we show the sufficiency part. Let $\det\big(b_{ij}\big)_{1\leq i, j\leq k}\neq0$.
Taking any two elements $\alpha, \beta \in {\bf F}_{q^m}$ such that
$F(\alpha)=F(\beta)$, i.e.,

$$L(\alpha)+\sum_{j=1}^k\gamma_jh_j(f_j(\alpha))
=L(\beta)+\sum_{j=1}^k\gamma_jh_j(f_j(\beta)),$$
we then derive that
$$L(\alpha-\beta)=\sum_{j=1}^k\gamma_j\big(h_j(f_j(\beta))-h_j(f_j(\alpha))\big).\eqno(2.1)$$
Since $\gamma_j\in {\rm Ker}(L)$ and $h_j(f_j(\beta))-h_j(f_j(\alpha))\in{\bf F}_{q}$
for $1\leq j\leq k$, one deduces that
$\sum_{j=1}^k\gamma_j\big(h_j(f_j(\beta))-h_j(f_j(\alpha))\big)\in {\rm Ker}(L)$. By (2.1),
we know that $\sum_{j=1}^k\gamma_j\big(h_j(f_j(\beta))-h_j(f_j(\alpha))\big)\in {\rm Im}(L)$.
But ${\rm Ker}(L)\cap {\rm Im}(L)=\{0\}$. So
$$\sum_{j=1}^k\gamma_j\big(h_j(f_j(\beta))-h_j(f_j(\alpha))\big)=0.\eqno(2.2)$$
By (2.1) and (2.2), we have $L(\alpha-\beta)=0$. So $\alpha-\beta\in{\rm Ker}(L)$.
Since $\{\gamma_1,...,\gamma_k\}$ is a basis of ${\rm Ker}(L)$ over ${\bf F}_q$,
it then follows that there exist $a_1,...,a_k\in {\bf F}_q$ such that
$$\alpha=\beta+a_1\gamma_1+...+a_k\gamma_k.\eqno(2.3)$$
Notice that $\{\gamma_1,...,\gamma_k\}$ is a basis of ${\rm Ker}(L)$ over ${\bf F}_q$,
we know that $\gamma_1,...,\gamma_k$ are linearly independent over ${\bf F}_q$.
Then by (2.2), we have for $1\leq j\leq k$ that
$$h_j(f_j(\beta))-h_j(f_j(\alpha))=0.\eqno(2.4)$$
Replacing $\alpha$ by $\beta+a_1\gamma_1+...+a_k\gamma_k$ in
(2.4) gives us that
$$h_j(f_j(\beta))-h_j(f_j(\beta+a_1\gamma_1+...+a_k\gamma_k))=0\ {\rm for}\ 1\leq
j\leq k. \eqno(2.5)$$
Since $h_j(x)$ is a permutation polynomial of
${\bf F}_q$, (2.5) is equivalent to
$$f_j(\beta+a_1\gamma_1+...+a_k\gamma_k)-f_j(\beta)=0\ {\rm for}\ 1\leq
j\leq k.\eqno (2.6)$$

On the other hand, since $\gamma_i$ is a $b_{ij}$-linear translator
of $f_j$ for all $1\le i, j\le k$, we can deduce that
$f_j(\beta+a_1\gamma_1+...+a_k\gamma_k)-f_j(\beta)=a_1b_{1j}+a_2b_{2j}...+a_kb_{kj}$.
Thus (2.6) is equivalent to
$$a_1b_{1j}+a_2b_{2j}...+a_kb_{kj}=0 \ {\rm for}\ 1\leq
j\leq k. \eqno(2.7)$$
It follows that $(a_1,...,a_k)\in {\bf F}_q^k$ is a solution
of the following system of linear equations:
\[\begin{cases}x_1b_{11}+x_2b_{21}+...+x_kb_{k1}=0\\
   x_1b_{12}+x_2b_{22}+...+x_kb_{k2}=0\\
   \ \ \ \ \ \ \ \ \ \ \ \ \ \ \ \vdots \\
   x_1b_{1k}+x_2b_{2k}+...+x_kb_{kk}=0.
   \end{cases}\eqno(2.8) \]
So by $\det\big(b_{ij}\big)_{1\leq i, j\leq k}\neq0$ we know that
the rank of the coefficient matrix of (2.8) is equal to $k$.
It follows that the system (2.8) of linear equations
has only zero solution. Namely, $(a_1, ..., a_k)=(0, ..., 0)$. So
by (2.3), we get that $\alpha=\beta$. Therefore $F(x)$ is a permutation
polynomial of ${\bf F}_{q^m}$. The sufficiency part is proved.

Let us now show the necessity part. Let $F(x)$ be a permutation
polynomial of ${\bf F}_{q^m}$. Suppose that $(a_1,...,a_k)\in
{\bf F}_q^k$ is a solution the system (2.8) of linear equations. Then
(2.7) is satisfied. By the equivalence of (2.5) and (2.7), we obtain that
$$h_j(f_j(\beta))-h_j(f_j(\beta+a_1\gamma_1+...+a_k\gamma_k))=0\ {\rm for}\ 1\leq j\leq k,$$
where $\beta\in{\bf F}_{q^m}$. Writing $\alpha:=\beta+a_1\gamma_1+...+a_k\gamma_k$
gives us that
$$\sum_{j=1}^k\gamma_j\big(h_j(f_j(\beta))-h_j(f_j(\alpha))\big)=0.$$
Since $\gamma_1, ..., \gamma_k\in {\rm Ker}(L)$, it follows that
$$L(\alpha-\beta)=L(a_1\gamma_1+...+a_k\gamma_k)=a_1L(\gamma_1)+...+a_kL(\gamma_k)=0.$$
It then follows that
$$L(\alpha-\beta)=\sum_{j=1}^k\gamma_j\big(h_j(f_j(\beta))-h_j(f_j(\alpha))\big).$$
We can derive immediately that $F(\alpha)=F(\beta)$.
Since $F(x)$ is a permutation polynomial of ${\bf F}_{q^m}$,
we have $\alpha=\beta$. Hence $a_1\gamma_1+...+a_k\gamma_k=0.$
But $\{\gamma_1,...,\gamma_k\}$ is a basis of ${\rm Ker}(L)$ over
${\bf F}_q$. Thus $a_1=...=a_k=0$. That is, the system (2.8)
of linear equations has only zero solution. Thus
$\det\big(b_{ij}\big)_{1\leq i, j\leq k}\neq0$
as desired. The necessity part is proved.

The proof of Theorem 2.1 is complete.
\end{proof}

By Theorem 2.1, we get the following interesting results.\\

\noindent{\bf Corollary 2.1.} {\it Let $m\geq2$ be a positive integer
with ${\rm gcd}(p, m)=1$, $\gamma_1,...,\gamma_{m-1}\in{\bf F}_{q^m}\setminus {\bf F}_{q}$
be linearly independent over ${\bf F}_q$ and $h_1(x), ..., h_{m-1}(x)\in {\bf F}_{q}[x]$
be permutation polynomials of ${\bf F}_{q}$. For any integers $i$ and $j$ with
$1\leq i, j\leq m-1$, let  $b_{ij}\in {\bf F}_q$ and $\gamma_i$ be a $b_{ij}$-linear
translator of $f_j:{\bf F}_{q^m}\to {\bf F}_{q}$.
Then $F(x):=\textup{Tr}_{{\bf F}_{q^m}/{\bf F}_{q}}(x)+\sum_{j=1}^{m-1}\gamma_jh_j(f_j(x))$
is a permutation polynomial of ${\bf F}_{q^m}$ if and only if $\det\big(b_{ij}\big)_{1\leq i,j\leq m-1}\neq0$.}
\begin{proof}
Since $\textup{Tr}_{{\bf F}_{q^m}/{\bf F}_{q}}:{\bf F}_{q^m}\to {\bf F}_{q}$ is surjective,
one has ${\rm Im}(\textup{Tr}_{{\bf F}_{q^m}/{\bf F}_{q}})={\bf F}_{q}$. For any
$u\in{\rm Ker}(\textup{Tr}_{{\bf F}_{q^m}/{\bf F}_{q}})\cap{\bf F}_{q}$,
we have $\textup{Tr}_{{\bf F}_{q^m}/{\bf F}_{q}}(u)=0$ and $\textup{Tr}_{{\bf F}_{q^m}/{\bf F}_{q}}(u)=mu$.
Thus $mu=0$. But the hypothesis that ${\rm gcd}(p, m)=1$ implies that $u=0$. Thus
${\rm Ker}(\textup{Tr}_{{\bf F}_{q^m}/{\bf F}_{q}})\cap{\bf F}_{q}=\{0\}.$ So
${\rm Ker}(\textup{Tr}_{{\bf F}_{q^m}/{\bf F}_{q}})={\bf F}_{q^m}\setminus {\bf F}_{q}^*$.
Then applying Theorem 2.1 to
$\textup{Tr}_{{\bf F}_{q^m}/{\bf F}_{q}}(x)$ concludes Corollary 2.1.
\end{proof}

\noindent{\bf Corollary 2.2.} {\it Let $p$ be an odd prime and $k$ be a positive integer.
Let $\{\gamma_1,...,\gamma_k\}$ be a basis of ${\bf F}_{q^k}$ over ${\bf F}_q$
and $h_1(x), ..., h_k(x)\in {\bf F}_{q}[x]$ be permutation polynomials of ${\bf F}_{q}$.
For any integers $i$ and $j$ with $1\leq i,j\leq k$, let  $b_{ij}\in {\bf F}_q$
and $\gamma_i$ be a $b_{ij}$-linear translator of $f_j:{\bf F}_{q^{2k}}\to {\bf F}_{q}$.
Then $F(x):=x-x^{q^k}+\sum_{j=1}^k\gamma_{j}h_j(f_j(x))$ is a
permutation polynomial of ${\bf F}_{q^{2k}}$ if and only if $\det\big(b_{ij}\big)_{1\leq i, j\leq k}\ne0$.}
\begin{proof}
For any $u\in{\rm Ker}(x-x^{q^k})\cap {\rm Im}(x-x^{q^k})$, we have $u=u^{q^k}$ and $u=v-v^{q^k}$
for some $v\in{\bf F}_{q^{2k}}$. It follows that
$u=(v-v^{q^k})^{q^k}=v^{q^k}-v^{q^{2k}}=v^{q^k}-v=-u$. It implies that
$2u=0$. But $p$ is an odd prime. So $u=0$. We conclude that
${\rm Ker}(x-x^{q^k})\cap {\rm Im}(x-x^{q^k})=\{0\}.$
So setting $L(x)=x-x^{q^k}$ and $m=2k$ in Theorem 2.1 gives us Corollary 2.2.
\end{proof}

By Lemma 2.1, we can construct some special mappings
$f:{\bf F}_{q^m}\to {\bf F}_{q}$ having linear translators. Thus
Corollaries 2.1 and 2.2 give us the following interesting consequences.\\

\noindent{\bf Corollary 2.3.} {\it Let $m\geq2$ be a positive integer
with ${\rm gcd}(p, m)=1$, $\gamma_1,...,\gamma_{m-1}\in{\bf F}_{q^m}\setminus{\bf F}_{q}$
be linearly independent over ${\bf F}_q$.
Let $h_j: {\bf F}_{q}\to {\bf F}_{q}$ be a permutation
of ${\bf F}_{q}$ and $H_j:{\bf F}_{q^m}\to {\bf F}_{q^m}, \beta_j\in {\bf F}_{q^m}$
for $1\leq j\leq m-1$. Then $F(x):=\textup{Tr}_{{\bf F}_{q^m}/{\bf F}_{q}}(x)
+\sum_{j=1}^{m-1}\gamma_jh_j(\textup{Tr}_{{\bf F}_{q^m}/{\bf
F}_{q}}(H_j(\textup{Tr}_{{\bf F}_{q^m}/{\bf F}_{q}}(x))+\beta_jx))$
is a permutation polynomial of ${\bf
F}_{q^m}$ if and only if $\det\big(\textup{Tr}_{{\bf F}_{q^m}/{\bf F}_{q}}
(\gamma_i\beta_j)\big)_{1\leq i,j\leq m-1}\ne0$.}

\begin{proof}
In Corollary 2.1, we set $f_j(x)=\textup{Tr}_{{\bf
F}_{q^{m}}/{\bf F}_{q}}(H_j(\textup{Tr}_
{{\bf F}_{q^m}/{\bf F}_{q}}(x))+\beta_jx)$. It is
easy to check that $\gamma_i$ is a
$\textup{Tr}_{{\bf F}_{q^m}/{\bf F}_{q}}(\gamma_i\beta_j)$-linear
translator of $f_j(x)$ for $1\le i, j \le m-1$. Then it follows immediately
from Corollary 2.1 that $F(x)$ is a permutation polynomial of ${\bf
F}_{q^m}$ if and only if
$\det\big(\textup{Tr}_{{\bf F}_{q^m}/{\bf F}_{q}}(\gamma_i\beta_j)
\big)_{1\leq i,j\leq m-1}\ne0$. Hence Corollary 2.3 is proved.
\end{proof}

\noindent{\bf Corollary 2.4.} {\it Let $p$ be an odd prime and $k$ be a positive integer.
Let $\alpha\in {\bf F}_{q^k}$ be a primitive element of ${\bf F}_{q^k}$.
Let $h_1(x), ..., h_k(x)\in {\bf F}_{q}[x]$ be permutation
polynomials of ${\bf F}_{q}$, $H_1(x),..., H_k(x)\in {\bf F}_{q^{2k}}[x]$ and
$\beta_1,...,\beta_k\in{\bf F}_{q^{2k}}$.
Then $F(x):=x-x^{q^k}+\sum_{j=1}^k\alpha^{j-1}h_j(\textup{Tr}_{{\bf F}_{q^{2k}}/{\bf
F}_{q}}(H_j(x-x^{q^k})+\beta_jx))$ is a
permutation polynomial of ${\bf F}_{q^{2k}}$ if and only if
$\det\big(\textup{Tr}_{{\bf F}_{q^{2k}}/{\bf F}_{q}}(\alpha^{i-1}\beta_j)\big)_{1\leq i, j\leq k}\ne0$.}
\begin{proof}
Since $\alpha\in {\bf F}_{q^k}$ is a primitive element of ${\bf F}_{q^k}$,
it follows that the set $\{1,\alpha,...,\alpha^{k-1}\}$ is a basis of ${\bf F}_{q^k}$.
It is easy to check that $\alpha^{i-1}$ is a $\textup{Tr}_{{\bf F}_{q^{2k}}/{\bf
F}_{q}}(\alpha^{i-1}\beta_j)$-linear translator of $\textup{Tr}_{{\bf F}_{q^{2k}}/{\bf
F}_{q}}(H_j(x-x^{q^k})+\beta_jx)$ for $1\le i,j\le k$.
Applying Corollary 2.2 to $f_j=\textup{Tr}_{{\bf F}_{q^{2k}}/{\bf
F}_{q}}(H_j(x-x^{q^k})+\beta_jx)$ and $\gamma_j=\alpha^{j-1}$ for
$1\le j\le k$ gives us that $F(x)$ is a permutation polynomial
of ${\bf F}_{q^{2k}}$ if and only if $\det\big(\textup{Tr}_{{\bf F}_{q^{2k}}/
{\bf F}_{q}}(\alpha^{i-1}\beta_j)\big)_{1\leq i, j\leq k}\ne0$.
\end{proof}
 To illustrate Corollaries 2.3 and 2.4, we give the following examples.\\

\noindent{\bf Example 2.1.} {\it Let $p$ be an odd prime and $t_1, t_2$ be
positive integers satisfying that $\textup{gcd}(t_i, q-1)=1$ for $i=1,2$.
Let $\alpha\in{\bf F}_{q^2}\setminus{\bf F}_q$, $\beta_1, \beta_2\in {\bf F}_{q^{4}}$
and $H_1(x), H_2(x)\in {\bf F}_{q^{4}}[x]$.
Then $F(x):=x^{q^2}-x+(\textup{Tr}_{{\bf F}_{q^{4}}/{\bf
F}_{q}}(H_1(x^{q^2}-x)+\beta_1x))^{t_1}+\alpha(\textup{Tr}_{{\bf
F}_{q^{4}}/{\bf F}_{q}}(H_2(x^{q^2}-x)+\beta_2x))^{t_2}$ is a
permutation polynomial of ${\bf F}_{q^{4}}$ if and only if
\begin{equation*}
\det\begin{pmatrix} \textup{Tr}_{{\bf F}_{q^{4}}/{\bf F}_{q}}(\beta_1) & \textup{Tr}_{{\bf F}_{q^{4}}
 /{\bf F}_{q}}(\beta_2) \\
 \textup{Tr}_{{\bf F}_{q^{4}}/{\bf
F}_{q}}(\alpha\beta_1) & \textup{Tr}_{{\bf F}_{q^{4}}/{\bf
F}_{q}}(\alpha\beta_2) \end{pmatrix}\ne 0.
\end{equation*}}

\noindent{\bf Example 2.2.} {\it Let $p$ be an odd prime and $t_1, t_2, t_3$ be
positive integers satisfying that $\textup{gcd}(t_i, q^2-1)=1$ for $i=1,2,3$. Let
$\beta_1, \beta_2, \beta_3\in {\bf F}_{q^{4}}$ and $H_1(x), H_2(x),H_3(x)\in {\bf F}_{q^{4}}[x]$.
Let $\alpha\in{\bf F}_{q^4}$ be a primitive element of ${\bf F}_{q^4}$ and $D_{t_i}(x,1)$
be a Dickson polynomial for $i=1,2,3$.
Then $F(x):=x^{q^3}+x^{q^2}+x^q+x+\sum_{i=1}^3\alpha^i D_{t_i}(\textup{Tr}_{{\bf F}_{q^{4}}/{\bf
F}_{q}}(H_i(x^{q^3}+x^{q^2}+x^q+x)+\beta_ix),1)$ is a
permutation polynomial of ${\bf F}_{q^{4}}$ if and only if
\begin{equation*}
\det\begin{pmatrix} \textup{Tr}_{{\bf F}_{q^{4}}/{\bf F}_{q}}(\alpha\beta_1) &
\textup{Tr}_{{\bf F}_{q^{4}}/{\bf F}_{q}}(\alpha\beta_2) &
\textup{Tr}_{{\bf F}_{q^{4}}/{\bf F}_{q}}(\alpha\beta_3)\\
\textup{Tr}_{{\bf F}_{q^{4}}/{\bf F}_{q}}(\alpha^2\beta_1)&
\textup{Tr}_{{\bf F}_{q^{4}}/{\bf F}_{q}}(\alpha^2\beta_2)&
\textup{Tr}_{{\bf F}_{q^{4}}/{\bf F}_{q}}(\alpha^2\beta_3)\\
\textup{Tr}_{{\bf F}_{q^{4}}/{\bf F}_{q}}(\alpha^{3}\beta_1) &
\textup{Tr}_{{\bf F}_{q^{4}}/{\bf F}_{q}}(\alpha^{3}\beta_2)&
\textup{Tr}_{{\bf F}_{q^{4}}/{\bf F}_{q}}(\alpha^{3}\beta_3) \end{pmatrix}\ne 0.
\end{equation*}}

We are now in a position to state the second main result of this paper.\\

\noindent{\bf Theorem 2.2.} {\it Let $k$ and $l$ be positive integers with $l\leq k$.
For any integers $i$ and $j$ with $1\le i,j \le k$, let $\gamma_i\in {\bf F}_{q^m}$,
$b_{ij}\in{\bf F}_{q}$ and $\gamma_i$  be a $b_{ij}$-linear translator of $f_j:{\bf F}
_{q^m}\to {\bf F}_q$ such that $\gamma_1,...,\gamma_k$ are linearly independent over ${\bf F}_q$.
Let $A=\big(b_{ij}\big)_{1\leq i, j\leq k}$ be a $k\times k$ matrix
over ${\bf F}_q$ and $I$ be the $k\times k$ identity
matrix over ${\bf F}_q$. Then each of the following is true:\\
{\rm (1)} $F(x):=x+\sum_{j=1}^k\gamma_jf_j(x)$ is a permutation
polynomial of ${\bf F}_{q^m}$ if and only if $\textup{rank}(I+A)=k$. \\
{\rm (2)} $F(x):=x+\sum_{j=1}^k\gamma_jf_j(x)$ is a
$q^l$-to-$1$ mapping of  ${\bf F}_{q^m}$ if
$\textup{rank}(I+A)=k-l$.}

\begin{proof}
(1) Assume that $\textup{rank}(I+A)=k$.
Take any two elements $\alpha, \beta \in {\bf F}_{q^m}$ satisfying
$F(\alpha)=F(\beta)$, that is,
$$\alpha+\sum_{j=1}^k\gamma_jf_j(\alpha)=\beta+\sum_{j=1}^k
\gamma_jf_j(\beta),\eqno(2.9)$$ which is equivalent to
$$\alpha-\beta=\sum_{j=1}^k\gamma_j\big(f_j(\beta)-f_j(\alpha)\big).\eqno(2.10)$$
Writing $a_j:=f_j(\beta)-f_j(\alpha)\in {\bf F}_q$, then by (2.10), we
get that $\alpha=\beta+\sum_{j=1}^k\gamma_ja_j.$ Replacing $\alpha$ by
$\beta+\sum_{j=1}^k\gamma_ja_j$ in (2.9), we arrive at
$$\sum_{j=1}^k\gamma_j\big(a_j+f_j(\beta+
\sum_{i=1}^k\gamma_ia_i)-f_j(\beta)\big)=0.\eqno(2.11)$$

Since $\gamma_i$ is a $b_{ij}$-linear translator of
$f_j$ for $1\le i,j \le k$, one has $f_j(\beta+
\sum_{i=1}^k\gamma_ia_i)-f_j(\beta)=\sum_{i=1}^k a_ib_{ij}$. Thus
(2.11) is equivalent to
$$\sum_{j=1}^k\gamma_j\big(a_j+
\sum_{i=1}^k a_ib_{ij}\big)=0.\eqno(2.12)$$ Since
$\gamma_1,...,\gamma_k$ are linearly independent over ${\bf F}_q$,
(2.12) is equivalent to
$$a_j+ \sum_{i=1}^k a_ib_{ij}=0 \ {\rm for} \ 1\le j\le k. \eqno(2.13)$$
Thus $(a_1, ..., a_k)^T\in {\bf F}_q^k$ is a solution of the system of linear equations
$$(I+A)^TX=0, \eqno(2.14)$$
where $(I+A)^T$ stands for the transpose of $I+A$ and $X=(x_1, ,..., x_k)^T$.

Since $\textup{rank}(I+A)=k$, the system (2.14) of linear equations
has only zero solution. Thus $a_1=...=a_k=0$. It follows from
$\alpha=\beta+\sum_{j=1}^k\gamma_ja_j$ that $\alpha=\beta$. Thus
$F(x)$ is a permutation polynomial of ${\bf F}_{q^m}$. So the
sufficiency part of (1) is proved.

Now we prove the necessity part of (1). Suppose that $F(x)$ is a
permutation polynomial of ${\bf F}_{q^m}$. If $(a_1, ..., a_k)^T\in {\bf F}_q^k$
is a solution of the system (2.14) of linear equations, then (2.13) is true. By
the equivalence between (2.13) and (2.11), we can deduce that
$$\beta+\sum_{j=1}^k\gamma_ja_j+\sum_{j=1}^k\gamma_jf_j
(\beta+\sum_{i=1}^k\gamma_ia_i)
=\beta+\sum_{j=1}^k\gamma_jf_j(\beta)$$ for $\beta\in{\bf F}_{q^m}$.
Putting $\alpha:=\beta+\sum_{j=1}^k\gamma_ja_j$ gives us that
$$\alpha+\sum_{j=1}^k\gamma_jf_j(\alpha)=\beta+\sum_{j=1}^k
\gamma_jf_j(\beta).$$
In other words, one has $F(\alpha)=F(\beta)$. Since $F(x)$ is a permutation
polynomial of ${\bf F}_{q^m}$, we have $\alpha=\beta$.
It implies that $\sum_{j=1}^k\gamma_ja_j=0$.
But $\gamma_1,...,\gamma_k$ are linearly independent over ${\bf F}_q$,
we have $(a_1, ..., a_k)^T=(0, ..., 0)^T$. Thus the system (2.14)
of linear equations has only zero solution. So
$\textup{rank}(I+A)=k$. The necessity part of (1) is proved.

(2) Let $\textup{rank}(I+A)=k-l$. If $(c_1,...,c_k)^T\in {\bf F}_q^k$ is any solution
of the system (2.14) of linear equations and $\beta\in{\bf F}_{q^m}$, then
$c_j+\sum_{i=1}^kc_ib_{ij}=0$. It follows that
\begin{align*}
&F(\beta+\sum_{j=1}^k\gamma_jc_j)\\
=&\beta+\sum_{j=1}^k\gamma_jc_j+ \sum_{j=1}^k\gamma_jf_j(\beta+\sum_{i=1}^k\gamma_ic_i)\\
=&\beta+\sum_{j=1}^k\gamma_jc_j+\sum_{j=1}^k\gamma_j(f_j(\beta)+\sum_{i=1}^k
c_ib_{ij})\ ({\rm since}\ \gamma_i\ {\rm is\ a }\ b_{ij}-{\rm linear\ translator\ of } f_j)\\
=&\beta+\sum_{j=1}^k\gamma_jf_j(\beta)+\sum_{j=1}^k\gamma_j(c_j+\sum_{i=1}^kc_ib_{ij})\\
=&\beta+\sum_{j=1}^k\gamma_jf_j(\beta)=F(\beta).
\end{align*}

On the other hand, since $\textup{rank}(I+A)=k-l$, we know that
the dimension of the space of the solutions
of the system (2.14) of linear equations over ${\bf F}_q$ equals $l$,
(2.14) has exactly $q^l$ solutions. Since $\gamma_1, ..., \gamma_k$
are linearly independent over ${\bf F}_q$, it follows that
$$\#\{\sum_{j=1}^k\gamma_jc_j: (c_1, ..., c_k)^T \ {\rm satisfies}
\ {\rm that} \ (2.14)\}=q^l.$$
Therefore $F(x)$ is a $q^l$-to-$1$ mapping of ${\bf F}_{q^m}$. So part (2) is proved.

This completes the proof of Theorem 2.2.\end{proof}

The referee pointed out that part (2) of Theorem 2.2 has appeared in Theorem 3
of \cite{[EKK]}. We note that there are two typos in the statement of Theorem 3
of \cite{[EKK]}. That is, ``$(a_1, ..., a_k)^T\in {\bf F}_q^n$" should read as
``$(a_1, ..., a_k)^T\in {\bf F}_q^k$", and ``the mapping $F$ is a $q^{n-r}$-to-1 on
${\bf F}_{q^n}$" should read as ``the mapping $F$ is $q^{k-r}$-to-1 on ${\bf F}_{q^n}$".
Now picking $k=1$ and $l=1$, we then have the following result due to Kyureghyan \cite{[K]}.\\

\noindent{\bf Corollary 2.5.} \cite{[K]} {\it Let $\gamma \in {\bf
F}_{q^m}$ be a $b$-linear translator of $f:{\bf F}_{q^m}\to {\bf
F}_q$. Then each of the following is true.\\
{\rm (1)} $F(x):=x+\gamma f(x)$ is a permutation polynomial of
${\bf F}_{q^m}$, if $b\ne -1$.\\
{\rm (2)} $F(x):=x+\gamma f(x$ is a $q$-to-$1$ mapping of
${\bf F}_{q^m}$, if $b=-1$.}\\

For $k=2$, Kyureghyan \cite{[K]} gave the following results.\\

\noindent{\bf Corollary 2.6.} \cite{[K]} {\it Let $\gamma, \delta\in
{\bf F}_{q^m}$ be linearly independent over ${\bf F}_q$. Suppose
$\gamma$ is a $b_1$-linear translator of $f:{\bf F}_{q^m}\to {\bf
F}_q$ and a $b_2$-linear translator of $g:{\bf F}_{q^m}\to {\bf
F}_q$ and moreover $\delta$ is a $d_1$-linear translator of $f:{\bf
F}_{q^m}\to {\bf F}_q$ and a $d_2$-linear translator of $g:{\bf
F}_{q^m}\to {\bf F}_q$. Then $F(x):=x+\gamma f(x)+\delta g(x)$ is a
permutation polynomial of ${\bf F}_{q^m}$, if $b_1\ne -1$ and
$d_2-\frac{d_1b_2}{1+b_1}\ne -1$ or by symmetry, if $d_2\ne -1$ and
$b_1-\frac{d_1b_2}{1+d_2}\ne -1$.}\\

\noindent{\bf Corollary 2.7.} \cite{[K]} {\it Let $\gamma\in {\bf
F}_{q^m}\setminus{\bf F}_q$ and
$M(x):=x^{q^2}-(1+(\gamma^q-\gamma)^{q-1})x^q+(\gamma^q-\gamma)^{q-1}x.$
Let $H_1, H_2:{\bf F}_{q^m}\to {\bf F}_{q^m}$ and $\beta_1, \beta_2\in {\bf F}_{q^m}$.
Then $F(x):=x+\textup{Tr}_{{\bf F}_{q^m}/{\bf
F}_{q}}(H_1(M(x))+\beta_1 x)+\gamma\textup{Tr}_{{\bf F}_{q^m}/{\bf
F}_{q}}(H_1(M(x))+\beta_2 x)$ is a permutation polynomial of ${\bf
F}_{q^m}$ if $(1+\textup{Tr}_{{\bf F}_{q^m}/{\bf
F}_{q}}(\beta_1))\\(1+\textup{Tr}_{{\bf F}_{q^m}/{\bf
F}_{q}}(\gamma\beta_2)) \ne \textup{Tr}_{{\bf F}_{q^m}/{\bf
F}_{q}}(\gamma\beta_1)\textup{Tr}_{{\bf F}_{q^m}/{\bf
F}_{q}}(\beta_2).$}\\

As a special case of Theorem 2.2, we have the following interesting results.\\

\noindent{\bf Corollary 2.8.} {\it Let $k$ be a positive integer.
Let $L:{\bf F}_{q^m}\to {\bf F}_{q^m}$ be a linearized polynomial
with kernel ${\rm Ker}(L)$ and $\{\theta_1,...,\theta_k\}$
be a basis of ${\rm Ker}(L)$ over ${\bf F}_q$.
Let $H_j:{\bf F}_{q^m}\to {\bf F}_{q^m}$ and
$\beta_j\in {\bf F}_{q^m}$ for $1\le j\le k$.
Then $F(x):=x+\sum_{j=1}^k\theta_j\textup{Tr}_{{\bf F}_{q^m}/{\bf
F}_{q}}(H_j(L(x))+\beta_jx)$
is a permutation polynomial of ${\bf F}_{q^m}$ if and only if
$\det\big(I+\big(\textup{Tr}_{{\bf F}_{q^m}/{\bf F}_{q}}
(\theta_i\beta_j)\big)\big)_{1\leq  i, j\leq k}\ne 0.$}

\begin{proof}
Since $\{\theta_1,...,\theta_k\}$ is a basis of
${\rm Ker}(L)$ over ${\bf F}_q$, it follows that
$\theta_1,...,\theta_k$ are linearly independent over ${\bf F}_q$.
It is easy to check that $\theta_i$ is a $\textup{Tr}_{{\bf
F}_{q^m}/{\bf F}_{q}}(\theta_i\beta_j)$-linear translator of
$\textup{Tr}_{{\bf F}_{q^m}/{\bf F}_{q}}(H_j(L(x))+\beta_jx)$ for $1\leq i, j\leq k$.
Then letting $\gamma_j=\theta_j$ and $f_j(x)=\textup{Tr}_{{\bf F}_{q^m}/{\bf F}_{q}}
(H_j(L(x))+\beta_jx)$ in Theorem 2.2, Corollary 2.8 follows immediately.
\end{proof}

\noindent{\bf Corollary 2.9.} {\it  Let $\alpha\in {\bf F}_{q^m}$
be a primitive element of ${\bf F}_{q^m}$ and $m>3$ be a integer.
Let
$$a=\frac{(\alpha-\alpha^{q^3})(\alpha^{q^2}-\alpha)^{q-1}}{\alpha^{q^2}-\alpha^q},
b=\frac{(\alpha^{q^3}-\alpha)(\alpha-\alpha^{q})^{q^2-1}}
{\alpha^{q^2}-\alpha^q}, c=-1-a-b,$$
and $N(x):=x^{q^3}+ax^{q^2}+bx^q+cx$. Let $H_1(x), H_2(x), H_3(x)\in
{\bf F}_{q^m}[x]$ and $\gamma_1, \gamma_2, \gamma_3\in{\bf F}_{q^m}$.
Then $F(x):=x+\textup{Tr}_{{\bf F}_{q^m}/{\bf F}_{q}}(H_1(N(x))+\gamma_1x)+\alpha
\textup{Tr}_{{\bf F}_{q^m}/{\bf F}_{q}}(H_2(N(x))+\gamma_2x)+\alpha^2
\textup{Tr}_{{\bf F}_{q^m}/{\bf F}_{q}}(H_3(N(x))+\gamma_3x)$
is a permutation polynomial of ${\bf F}_{q^m}$ if and only if
\begin{equation*}
\det\begin{pmatrix}
    1+\textup{Tr}_{{\bf F}_{q^m}/{\bf
F}_{q}}(\gamma_1) & \textup{Tr}_{{\bf F}_{q^m}/{\bf
F}_{q}}(\gamma_2)&  \textup{Tr}_{{\bf F}_{q^m}/{\bf
F}_{q}}(\gamma_3)\\
    \textup{Tr}_{{\bf F}_{q^m}/{\bf
F}_{q}}(\alpha\gamma_1) &1+\textup{Tr}_{{\bf F}_{q^m}/{\bf
F}_{q}}(\alpha\gamma_2) &\textup{Tr}_{{\bf F}_{q^m}/{\bf
F}_{q}}(\alpha\gamma_2) \\
\textup{Tr}_{{\bf F}_{q^m}/{\bf F}_{q}}(\alpha^2\gamma_1)
&\textup{Tr}_{{\bf F}_{q^m}/{\bf F}_{q}}(\alpha^2\gamma_2)
&1+\textup{Tr}_{{\bf F}_{q^m}/{\bf
F}_{q}}(\alpha^2\gamma_2) \\
\end{pmatrix}\ne 0.
\end{equation*}}

\begin{proof}
Since $\alpha\in {\bf F}_{q^m}$ is a primitive element of ${\bf
F}_{q^m}$, $1,\alpha,\alpha^2$ are linearly independent over ${\bf
F}_q$. One can easily check that $1,\alpha,\alpha^2$ are the roots
of  $N(x)$ and $\alpha^i$ is a $\textup{Tr}_{{\bf F}_{q^m}/{\bf
F}_{q}}(\alpha^i\gamma_j)$ -linear translator of $\textup{Tr}_{{\bf
F}_{q^m}/{\bf F}_{q}}(H_j(N(x))+\gamma_jx)$ for $0\leq i\leq 2$ and
$1\le j\le 3$. Thus Corollary 2.9 follows immediately from Theorem 2.2.
\end{proof}

As an application of Theorem 2.2, we can get a large family
of {\it complete mappings} (also called {\it complete permutation 
polynomials}), which are the permutation polynomials $F(x)$
with $F(x)+x$ being a permutation polynomial as well.\\

\noindent{\bf Corollary 2.10.} {\it Let $p$ be an odd prime and $k$
be a positive integer. For any integers
$i$ and $j$ with $1\le i,j \le k$, let $\gamma_i\in {\bf F}_{q^m}$,
$b_{ij}\in{\bf F}_{q}$, $\gamma_i$  be a $b_{ij}$-linear translator
 of $f_j:{\bf F}_{q^m}\to {\bf F}_q$ such that
$\gamma_1,...,\gamma_k$ are linearly independent over ${\bf F}_q$.
Let $A=\big(b_{ij}\big)_{1\leq i, j\leq k}$ be a
$k\times k$ matrix over ${\bf F}_q$ and $I$ be the $k\times k$ identity matrix
over ${\bf F}_q$. Then $F(x):=x+\sum_{j=1}^k\gamma_jf_j(x)$ is a
complete mapping of ${\bf F}_{q^m}$ if and only if
$\textup{rank}(I+A)=k$ and $\textup{rank}(2I+A)=k$.}
\begin{proof}
In the similar way as in the proof of Theorem 2.2, we can show that
$2x+\sum_{j=1}^k\gamma_jf_j(x)$ is a permutation polynomial of ${\bf
F}_{q^m}$ if and only if $\textup{rank}(2I+A)=k$. Thus
$F(x):=x+\sum_{j=1}^k\gamma_jf_j(x)$ is a complete mapping of ${\bf
F}_{q^m}$ if and only if $\textup{rank}(I+A)=k$ and
$\textup{rank}(2I+A)=k$. So Corollary 2.10 is proved.
\end{proof}

\section{\bf Permutation polynomials of the form
$L(x)+\sum_{i=1}^l\gamma_i \textup{Tr}_{{\bf F}_{q^m}/{\bf F}_{q}}(h_i(x))$}

In \cite{[CK]}, Charpin and Kyureghyan studied permutation
polynomials of the type $F(x):=G(x)+\gamma \textup{Tr}_{{\bf
F}_{q}/{\bf F}_{p}}(h(x))$. When $G(x)$ is a permutation polynomial
or a linearized polynomial, they characterized permutation
polynomials of this shape. In this section, we characterize
permutation polynomials of the form $L(x)+\sum_{i=1}^l\gamma_i
\textup{Tr}_{{\bf F}_{q^m}/{\bf F}_{q}}(h_i(x))$.
The third main result of this paper is given as follows.\\

\noindent{\bf Theorem 3.1.} {\it  Let $l$ and $k$ be positive integers with $l\le k$.
Let $L(x)\in{\bf F}_{q^m}[x]$ be a linearized polynomial such that $\dim({\rm Ker}(L))=k$ and
${\rm Ker}(L)\cap {\rm Im}(L)=\{0\}$. Let $\gamma_1,...,\gamma_l\in {\rm Ker}(L)$ be linearly
independent over ${\bf F}_q$ and $h_1(x), ..., h_l(x)\in {\bf F}_{q^m}[x]$.
Then $F(x):=L(x)+\sum_{i=1}^l\gamma_i \textup{Tr}_{{\bf F}_{q^m}/{\bf F}_{q}}(h_i(x))$ is a
permutation polynomial of ${\bf F}_{q^m}$ if and only if there exists an integer
$i$ with $1\le i \le l$ such that $\textup{Tr}_{{\bf F}_{q^m}/{\bf F}_{q}}
(h_i(x+\varepsilon)-h_i(x))\ne 0$ for any $x\in {\bf F}_{q^m}$ and any
$\varepsilon\in{\rm Ker}(L)\setminus\{0\}$.}

\begin{proof}
First we show the sufficiency part. Assume that there exists an integer
$i$ with $1\le i \le l$ such that
$\textup{Tr}_{{\bf F}_{q^m}/{\bf F}_{q}}(h_i(x+\varepsilon)-h_i(x))\ne 0$
for any $x\in {\bf F}_{q^m}$ and any $\varepsilon\in{\rm Ker}(L)\setminus\{0\}$.
Take any two elements $\alpha, \beta \in {\bf F}_{q^m}$ satisfying
$F(\alpha)=F(\beta)$, namely,
$$L(\alpha)+\sum_{i=1}^l\gamma_i \textup{Tr}_{{\bf F}_{q^m}/{\bf F}_{q}}(h_i(\alpha))=
L(\beta)+\sum_{i=1}^l\gamma_i \textup{Tr}_{{\bf F}_{q^m}/{\bf F}_{q}}(h_i(\beta)).$$
We deduce that
$$L(\alpha-\beta)=\sum_{i=1}^l\gamma_i\textup{Tr}_{{\bf F}_{q^m}/{\bf F}_{q}}
(h_i(\beta)-h_i(\alpha)).\eqno(3.1)$$
Since $\gamma_i\in {\rm Ker}(L)$ and $\textup{Tr}_{{\bf F}_{q^m}/{\bf F}_{q}}
(h_i(\beta)-h_i(\alpha))\in{\bf F}_{q}$ for $1\leq i\leq l$,
we get immediately that $\sum_{i=1}^l\gamma_i\textup{Tr}_{{\bf F}_{q^m}/{\bf F}_{q}}
(h_i(\beta)-h_i(\alpha))\in {\rm Ker}(L)$. But by (3.1), one has
$$\sum_{i=1}^l\gamma_i\textup{Tr}_{{\bf F}_{q^m}/{\bf F}_{q}}
(h_i\\(\beta)-h_i(\alpha))=L(\alpha-\beta)\in {\rm Im}(L).$$
It then follows from ${\rm Ker}(L)\cap {\rm Im}(L)=\{0\}$ that
$$\sum_{i=1}^l\gamma_i\textup{Tr}_{{\bf F}_{q^m}/{\bf F}_{q}}
(h_i(\beta)-h_i(\alpha))=0.\eqno(3.2)$$
Hence (3.1) together with (3.2) infers that $\alpha-\beta\in{\rm Ker}(L)$.
Thus there exists an element
$\varepsilon\in{\rm Ker}(L)$ such that $\alpha=\beta+\varepsilon$.

We claim that $\varepsilon=0$. Suppose that $\varepsilon\ne0$. By the hypothesis, we know that
there exists an integer $i_0$ with $1\le i_0 \le l$ such that $\textup{Tr}_{{\bf F}_{q^m}/{\bf F}_{q}}
(h_{i_0}(\beta+\varepsilon)-h_{i_0}(\beta))\ne 0$. Since $\gamma_1,...,\gamma_l$
are linearly independent over ${\bf F}_q$, it follows from (3.2) that for all $j$ with
$1\le j \le l$, one has $\textup{Tr}_{{\bf F}
_{q^m}/{\bf F}_{q}}(h_j(\beta)-h_j(\alpha))=0$, i.e.,
$\textup{Tr}_{{\bf F}_{q^m}/{\bf F}_{q}}(h_j(\beta+\varepsilon)-h_j(\beta))=0$.
In particular, we have
$\textup{Tr}_{{\bf F}_{q^m}/{\bf F}_{q}}(h_{i_0}(\beta+\varepsilon)-h_{i_0}(\beta))\ne 0$.
This arrives at a contradiction. Thus $\varepsilon=0$. The claim is proved.
Therefore $F(x)$ is a permutation polynomial of ${\bf F}_{q^m}$. The sufficiency part is proved.

Let us now show the necessity part. Let $F(x)$ be a permutation
polynomial of ${\bf F}_{q^m}$. For any $x\in{\bf F}_{q^m}$ and any
$\varepsilon\in{\rm Ker}(L)\setminus\{0\}$,
we have
$$F(x+\varepsilon)-F(x)=\sum_{i=1}^l\gamma_i \textup{Tr}_
{{\bf F}_{q^m}/{\bf F}_{q}}(h_i(x+\varepsilon)-h_i(x)).\eqno(3.3)$$
Since $F(x)$ is a permutation polynomial of ${\bf F}_{q^m}$, it follows from (3.3) that
$$\sum_{i=1}^l\gamma_i \textup{Tr}_{{\bf F}_{q^m}/{\bf F}_{q}}(h_i(x+\varepsilon)-h_i(x))\ne0.$$
But $\gamma_1,...,\gamma_l\in {\rm Ker}(L)$ are linearly independent over ${\bf F}_q$.
Hence there exists an integer $i$ with $1\le i \le l$ such that
$\textup{Tr}_{{\bf F}_{q^m}/{\bf F}_{q}}(h_i(x+\varepsilon)-h_i(x))\ne 0$ for any
$x\in {\bf F}_{q^m}$ and any $\varepsilon\in{\rm Ker}(L)\setminus\{0\}$.
The necessity part is proved.

The proof of Theorem 3.1 is complete.
\end{proof}

By Theorem 3.1, we can easily deduce the following consequences.\\

\noindent{\bf Corollary 3.1.} {\it Let $l$ and $m\geq2$ be positive integers
with ${\rm gcd}(p, m)=1$ and $l<m$. Let $\gamma_1,...,\gamma_{l}\in{\bf F}_{q^m}\setminus {\bf F}_{q}$
be linearly independent over ${\bf F}_q$ and $h_1(x), ..., h_{l}(x)\in {\bf F}_{q^m}[x]$.
Then $F(x):=\textup{Tr}_{{\bf F}_{q^m}/{\bf F}_{q}}(x)+\sum_{i=1}^l\gamma_i \textup{Tr}_
{{\bf F}_{q^m}/{\bf F}_{q}}(h_i(x))$ is a permutation polynomial of ${\bf F}_{q^m}$
if and only if there exists an integer $i$ with $1\le i \le l$ such that
$\textup{Tr}_{{\bf F}_{q^m}/{\bf F}_{q}}(h_i(x+\varepsilon)-h_i(x))\ne 0$
for any $x\in {\bf F}_{q^m}$ and any $\varepsilon\in{\bf F}_{q^m}\setminus {\bf F}_{q}$.}
\begin{proof}
By the proof of Corollary 2.1, we know that ${\rm Im}(\textup{Tr}_{{\bf F}_{q^m}/{\bf F}_{q}})={\bf F}_{q}$
 and ${\rm Ker}(\textup{Tr}_{{\bf F}_{q^m}/{\bf F}_{q}})\\\cap{\bf F}_{q}=\{0\}.$
Then applying Theorem 3.1 to
$\textup{Tr}_{{\bf F}_{q^m}/{\bf F}_{q}}(x)$ gives us Corollary 3.1.
\end{proof}

\noindent{\bf Corollary 3.2.} {\it Let $p$ be an odd prime, $l$ and $k$ be positive integers
with $l\le k$. Let $\{\gamma_1,...,\gamma_k\}$ be a basis of ${\bf F}_{q^k}$ over ${\bf F}_q$
and $h_1(x), ..., h_l(x)\in {\bf F}_{q^{2k}}[x]$.
Then $F(x):=x-x^{q^k}+\sum_{i=1}^l\gamma_i \textup{Tr}_{{\bf F}_{q^{2k}}/{\bf F}_{q}}(h_i(x))$ is a
permutation polynomial of ${\bf F}_{q^{2k}}$ if and only if there exists an integer $i$
with $1\le i \le l$ such that $\textup{Tr}_{{\bf F}_{q^{2k}}/{\bf F}_{q}}(h_i(x+\varepsilon)-h_i(x))\ne 0$
for any $x\in {\bf F}_{q^{2k}}$ and any $\varepsilon\in{\bf F}_{q^k}^*$.}
\begin{proof}
By the proof of Corollary 2.2, we conclude that ${\rm Ker}(x-x^{q^k})\cap {\rm Im}(x-x^{q^k})=\{0\}.$
So Corollary 3.2 follows from Theorem 3.1 by setting $L(x)=x-x^{q^k}$ and $m=2k$.
\end{proof}

{\bf Acknowledgement.} The authors thank the anonymous referee
for his/her careful reading of the manuscript and particularly for drawing
our attention to reference \cite{[EKK]}.

\end{document}